\documentclass{article}
\usepackage{amsmath}
\usepackage{amssymb}
\usepackage{amsthm}
\usepackage{amsfonts}

\numberwithin{equation}{section}

\theoremstyle{definition}
\newtheorem{theorem}{Theorem}[section]

\newtheorem*{thm}{Theorem}

\theoremstyle{definition}

\theoremstyle{remark}
\newtheorem{remark}{Remark}[section]

\begin{document}

\title{The Shatashvili-Vafa $G_{2}$ superconformal algebra as a Quantum Hamiltonian Reduction of $D(2,1;\alpha)$.}

\author{Reimundo Heluani\thanks{IMPA, Rio de Janeiro. Partially supported by CNPq
},\ L\'azaro O. Rodr\'iguez D\'iaz\thanks{IMPA, Rio de Janeiro. Supported by FAPERJ and CNPq.}.}

\date{}

\maketitle

\begin{abstract}
We obtain the superconformal algebra associated to a sigma model with target a manifold with $G_{2}$ holonomy, i.e., the Shatashvili-Vafa $G_{2}$ algebra as a quantum Hamiltonian reduction of the exceptional Lie superalgebra $D(2,1;\alpha)$ for $\alpha=1$. We produce the complete family of $W$-algebras $SW(\frac{3}{2},\frac{3}{2}, 2)$ (extensions of the $N=1$ superconformal algebra by two primary supercurrents of  conformal weight $\frac{3}{2}$ and $2$ respectively) as a quantum Hamiltonian reduction of $D(2,1;\alpha)$. As a corollary we find a free field realization of the Shatashvili-Vafa $G_{2}$ algebra, and an explicit description of the screening operators.
\end{abstract}

\section{Introduction}

The Shatashvili-Vafa $G_{2}$ algebra  \cite{Shatashvili-Vafa95} is a superconformal vertex algebra with six generators $\{L,G,\Phi,K,X,M\}$. It is an extension of the $N=1$ superconformal algebra of central charge $c=21/2$ (formed by the super-partners $\{L,G\}$) by two fields $\Phi$ and $K$, primary of conformal weight $\frac{3}{2}$ and $2$ respectively, and their superpartners $X$ and $M$ (of conformal weight $2$ and $\frac{5}{2}$ respectively). Their OPEs can be found in Appendix \ref{App:AppendixA} in the language of lambda brackets of \cite{dandrea}. 

This superconformal algebra appeared as the chiral algebra associated to the sigma model with target a manifold with $G_{2}$ holonomy in \cite{Shatashvili-Vafa95} its classical counterpart had been studied by Howe and Papadopoulos in \cite{Howe-Papadopoulos93}. In fact this algebra is a member of a two-parameter family $SW(\frac{3}{2}, \frac{3}{2}, 2)$ previously studied in \cite{Blumenhagen92} where the author found the family of all superconformal algebras which are extension of the super-Virasoro algebra, i.e., the $N=1$ superconformal algebra, by two primary supercurrents of conformal weights $\frac{3}{2}$ and $2$ respectively. It is a family  parametrized by $(c,\varepsilon)$ ($c$ is the central charge and $\varepsilon$ the coupling constant) of non-linear $W$-algebras. Its generators and relations are recalled in Appendix \ref{App:AppendixB}.
 
The Shatashvili-Vafa $G_{2}$ algebra is a quotient of $SW(\frac{3}{2},\frac{3}{2}, 2)$ with $c=\frac{21}{2}$ and $\varepsilon=0$, in other words is the only one among this family which has central charge $c=\frac{21}{2}$ and contains the tri-critical Ising model as a subalgebra.
It is precisely the fact that the Shatashvili-Vafa $G_{2}$ algebra appears as a $W$-algebra that motivated the authors to try to obtain this algebra as a quantum Hamiltonian reduction of some Lie superalgebra using the method developed in \cite{KacRoanWakimoto03}.

That $D(2,1;\alpha)$ is the right Lie superalgebra candidate to be used in the Hamiltonian reduction is known from scattered results in the physics literature. It was shown in \cite{Mallwitz95} that $SW(\frac{3}{2},\frac{3}{2}, 2)$ is the symmetry algebra of the quantized Toda theory corresponding to $D(2,1;\alpha)$ (in \cite{Nohara-Mohri91} was worked a classical version of this result in the case $\alpha=1$ ($D(2,1;\alpha)=osp(4|2)$)) and from the well established connection between the theory of nonlinear integrable equations and W-algebras, see for example \cite{FeiginFrenkel95}.

A coset realization of the $SW(\frac{3}{2},\frac{3}{2},2)$ superconformal algebra and therefore of the Shatashvili-Vafa algebra can be found in \cite{Noyvert02}. In \cite{FeiginSemikhatov01} was shown that the Hamiltonian reduction of $D(2,1;\alpha)$ coincides with this coset model (the authors however restrict their attention to the even part of the superalgebra). 

Some representations of the Shatashvili-Vafa $G_{2}$ superconformal algebra can be found in \cite{Noyvert02}, but the character formulae remains unknown. It was  observed in \cite{BoerNaqviShomer} that in order to systematically study the representation theory and the character formula for this algebra one should construct the Shatasvili-Vafa algebra using the quantum Drinfeld-Sokolov reduction developed in \cite{KacWakimoto04,FrenkelKacWakimoto92}. This step is accomplished in this article.

In section \ref{preliminaries} we review how to perform the quantum Hamiltonian reduction of a Lie superalgebra as introduced in \cite{KacRoanWakimoto03}. We recap some of the main theorems as well as under which conditions this Hamiltonian reduction process induces a free field realization.

In section \ref{quantum reduction of the algebra} we prove that the $SW(\frac{3}{2},\frac{3}{2},2)$ superconformal algebra is the quantum Hamiltonian reduction of the Lie superalgebra $D(2,1;\alpha)$, and obtain a free field realization of the $SW(\frac{3}{2},\frac{3}{2},2)$ algebra on a space of three free Bosons and three free Fermions. As particular cases $\alpha\in\{1,-\frac{1}{2},-2\}$ we obtain the Shatashvili-Vafa $G_{2}$ algebra as a quantum Hamiltonian reduction of the Lie superalgebra $osp(4|2)$, and also the corresponding free field realizations. We summarize our main result as (see Theorem \ref{thm:main-thm} its remark)

\noindent
\begin{thm}
 Let $\mathfrak{h}$ be the Cartan subalgebra of $D(2,1;\alpha)$. It is a three dimensional vector space with a non-degenerate bilinear form $(,)$ given by the Cartan matrix. Consider $\Pi \mathfrak{h}^*$ the odd vector space ($\Pi$ denotes parity change) with its natural bilinear form $-(,)$.  Let $V_{k}(\mathfrak{h}_{\mathrm{super}})$ be the super affine vertex algebra generated by three Bosons from $\mathfrak{h}$ and three Fermions from $\Pi \mathfrak{h^*}$ and lambda brackets
\[ [h_\lambda h'] = k \lambda (h,h'), \qquad [\phi_\lambda \phi']= -(\phi, \phi'), \qquad h,h' \in \mathfrak{h}, \: \phi, \phi' \in \Pi \mathfrak{h}^*. \]
\begin{enumerate}
\item $SW(\frac{3}{2},\frac{3}{2},2)$ is a sub-vertex algebra of $V_{k}(\mathfrak{h}_{\mathrm{super}})$. The particular case of central charge $c=21/2$ and vanishing coupling constant corresponds to $k=-2/3$. In this case the fields $G$ and $\Phi$ are given by \eqref{expressionforPhi} and all other fields can be obtained from these two. 
\item The Shatashvili-Vafa $G_2$ superconformal algebra is a quotient of this algebra by an ideal generated in conformal weight $7/2$ \eqref{ideal}. 
\item For each $\alpha_i$ of the three odd simple roots of $D(2,1, \alpha)$ there exists a module $M_i$ of $V_{k}(\mathfrak{h}_{\mathrm{super}})$ generated by a vector $|\alpha_i\rangle$ such that $h_n |\alpha_i\rangle = 0$ for $n > 0$, $h_0 |\alpha_i\rangle = (h, \alpha_i) |\alpha_i\rangle$ and $\phi_n |\alpha_i\rangle = 0$ for $n > 0$ (for all $h \in \mathfrak{h}$ and $\phi \in \Pi \mathfrak{h}^*$). Let $\Gamma_i(z)$ be the unique intertwiner of type $\binom{V_{k}(\mathfrak{h}_{\mathrm{super}})}{V_{k}(\mathfrak{h}_{\mathrm{super}}) \; M_i}$ and $Q_i \in \mathrm{Hom} (V_{k}(\mathfrak{h}_{\mathrm{super}}), M_i)$ its zero mode. Then for generic values of $(c,\varepsilon)$ we have $SW(\frac{3}{2},\frac{3}{2},2) = \cap_i Q_i \subset V_{k}(\mathfrak{h}_{\mathrm{super}})$.
\end{enumerate}
\end{thm}
In Section \ref{quantum reduction of the algebra} the reader can find a stronger version of this Theorem as the generators for $SW(\frac{3}{2},\frac{3}{2},2)$ are found for any values of the parameters $(c,\varepsilon)$.

\section{Quantum reduction of Lie superalgebras}\label{preliminaries}
In this section we recall the construction of the W-algebras $W_{k}(\mathfrak{g},x,f)$ introduced in \cite{KacRoanWakimoto03}. We follow the presentation in \cite{KacWakimoto04}. 

To construct the vertex algebra $W_{k}(\mathfrak{g},x,f)$ we need a quadruple $(\mathfrak{g},x,f,k)$ where $\mathfrak{g}=\mathfrak{g}_{\bar{0}}\oplus \mathfrak{g}_{\bar{1}}$ is a simple finite-dimensional Lie superalgebra with a non-degenerate even invariant supersymmetric bilinear form  $(.|.)$, and $x,f\in \mathfrak{g}_{\bar{0}}$ such that $ad\;x$ is diagonalizable on $\mathfrak{g}$ with half-integer eigenvalues, $[x,f]=-f$, the eigenvalues of ad $x$ on the centralizer $\mathfrak{g}^{f}$ of $f$ in $\mathfrak{g}$ are non-positive, and $k\in\mathbb{C}$.

We recall that a bilinear form $(.|.)$ on $\mathfrak{g}$ is called even if $(\mathfrak{g}_{\bar{0}}|\mathfrak{g}_{\bar{1}})=0$, supersymmetric if $(.|.)$ is symmetric (resp.\!\!\! skewsymmetric) on $\mathfrak{g}_{\bar{0}}$ (resp.\!\!\! $\mathfrak{g}_{\bar{1}}$), invariant if $([a,b]|c)=(a|[b,c])$ for all $a,b,c\in \mathfrak{g}$. 

A pair $(x,f)$ satisfying the above properties can be obtained when $x,f$ are part of an $\mathfrak{sl}_{2}$ triple, i.e., $[x,e]=e$, $[x,f]=-f$ and $[e,f]=x$. As this will be the case in the quantum reduction performed in section \ref{quantum reduction of the algebra}, we assume for the rest of this section that we are working with such a pair. Let $\mathfrak{g}=\oplus_{j\in\frac{1}{2}\mathbb{Z}}\mathfrak{g}_{j}$, be the eigenspace decomposition with respect to ad $x$. Denote 

\begin{equation}
\mathfrak{g}_{+}=\bigoplus_{j>0}\mathfrak{g}_{j},\;\;\;\; \mathfrak{g}_{-}=\bigoplus_{j<0}\mathfrak{g}_{j}, \;\;\;\; \mathfrak{g}_{\le}=\mathfrak{g}_{0}\bigoplus \mathfrak{g}_{-}.\nonumber
\end{equation}

Let $V_{k}(\mathfrak{g})$ denote the affine vertex algebra of level $k$ associated to $\mathfrak{g}$. Denote by $F(A)$ the vertex algebra of free superfermions associated to a vector superspace $A$ with an even skew-supersymmetric non-degenerate bilinear form $\left<.|.\right>$, i.e., the $\lambda$-bracket is given by $[\varphi_{\lambda}\psi]=\left<\varphi|\psi\right>$, $\varphi, \psi\in A$. 

\sloppy

On the vector superspace $\mathfrak{g}_{1/2}$ the element $f$ defines an even skew-supersymmetric non-degenerate bilinear form $\left<.|.\right>_{ne}$ by the formula:

\fussy

\begin{equation}
\left<a|b\right>=(f|[a,b]).\nonumber
\end{equation} 

The associated vertex algebra $F(\mathfrak{g}_{1/2})$ is called the vertex algebra of neutral free superfermions. Similary on the vector superspace $\Pi\mathfrak{g}_{+}\oplus\Pi\mathfrak{g}_{+}^{*}$( where $\Pi$ denotes parity-reversing), define an even skew-supersymmetric non-degenerate bilinear form $\left<.|.\right>_{ch}$ by:

\begin{equation}
\left<\Pi\mathfrak{g}_{+}|\Pi\mathfrak{g}_{+}\right>_{ch}=0=\left<\Pi\mathfrak{g}_{+}^*|\Pi\mathfrak{g}_{+}^*\right>_{ch},\nonumber
\end{equation}

\begin{equation}
\left<a|b^*\right>_{ch}=-(-1)^{p(a)p(b^*)}\left<b^*|a\right>_{ch}=b^*(a), \;\;\; a\in\Pi\mathfrak{g}_{+}, b^*\in \Pi\mathfrak{g}_{+}^*, \nonumber 
\end{equation}

where $p(a)$ denotes the parity of the element $a$. The associated vertex algebra $F(\Pi\mathfrak{g}_{+}\oplus\Pi\mathfrak{g}_{+}^{*})$ is called the vertex algebra of charged free superfermions. This vertex algebra carries an extra $\mathbb{Z}$-grading by charge by assigning: charge $\varphi=1$ and charge $\varphi^*=-1$, $\varphi\in\Pi\mathfrak{g}_{+}$, $\varphi^*\in\Pi\mathfrak{g}_{+}^{*}$. Consider the vertex algebra

\begin{equation}
C(\mathfrak{g},x,f,k)=V_{k}(\mathfrak{g})\otimes F(\Pi\mathfrak{g}_{+}\oplus\Pi\mathfrak{g}_{+}^{*})\otimes F(\mathfrak{g}_{1/2}). \nonumber
\end{equation} 

The charge decomposition of $F(\Pi\mathfrak{g}_{+}\oplus\Pi\mathfrak{g}_{+}^{*})$ induces a charge decomposition on $C(\mathfrak{g},x,f,k)$ by declaring charge $V_{k}(\mathfrak{g})=0$ and charge $F(\mathfrak{g}_{1/2})=0$. This makes $C(\mathfrak{g},x,f,k)$ a $\mathbb{Z}$-graded vertex algebra. We introduce a differential $d_{(0)}$ that makes $(C(\mathfrak{g},x,f,k),d_{(0)})$ a $\mathbb{Z}$-graded complex as follows. Let $\{u_{\alpha}\}_{\alpha\in S_{j}}$ be a basis of each $\mathfrak{g}_{j}$, an let $S:=\coprod_{j\in\frac{1}{2}\mathbb{Z}}S_{j}$, $S_{+}=\coprod_{j>0}S_{j}$. Put $m_{\alpha}:=j$ if $\alpha\in S_{j}$. The structure constants  $c_{\alpha\beta}^{\gamma}$ are defined by $[u_{\alpha},u_{\beta}]=\sum_{\gamma}c_{\alpha\beta}^{\gamma}u_{\gamma}$ for $(\alpha,\beta,\gamma\in S)$. Denote by $\{\varphi_{\alpha}\}_{\alpha\in S_{+}}$ the corresponding basis of $\Pi \mathfrak{g}_{+}$ and by $\{\varphi^{\alpha}\}_{\alpha\in S_{+}}$ the basis of $\Pi \mathfrak{g}_{+}^*$ such that $\left<\varphi_{\alpha},\varphi^{\beta}\right>_{ch}=\delta^{\beta}_{\alpha}$. Similary denote by $\{\Phi_{\alpha}\}_{\alpha\in S_{1/2}}$ the corresponding basis of $\mathfrak{g}_{1/2}$, and by $\{\Phi^{\alpha}\}_{\alpha\in S_{1/2}}$ the dual basis with respect to $\left<.,.\right>_{ne}$, i.e., $\left<\Phi_{\alpha},\Phi^{\beta}\right>_{ne}=\delta^{\beta}_{\alpha}$. It is useful to define $\Phi_{u}$ for any $u=\sum_{\alpha\in S}c_{\alpha}u_{\alpha}\in \mathfrak{g}$ by $\Phi_{u}:=\sum_{\alpha\in S_{1/2}}c_{\alpha}\Phi_{\alpha}$. Define the odd field
\begin{eqnarray}
d&=&\sum_{\alpha\in S_{+}}(-1)^{p(u_{\alpha})}:u_{\alpha}\varphi^{\alpha}:-\frac{1}{2}\sum_{\alpha,\beta,\gamma\in S_{+}}(-1)^{p(u_{\alpha})p(u_{\gamma})}c^{\gamma}_{\alpha \beta}:\varphi_{\gamma}\varphi^{\alpha}\varphi^{\beta}:
 \nonumber\\
&&+\sum_{\alpha\in S_{+}}(f|u_{\alpha})\varphi^{\alpha}+\sum_{\alpha\in S_{1/2}}:\varphi^{\alpha}\Phi_{\alpha}:.\nonumber
\end{eqnarray}

Its Fourier mode $d_{(0)}$ is an odd derivation of all products of the vertex algebra $C(\mathfrak{g},x,f,k)$, such that $d_{(0)}^2=0$ and that $d_{(0)}$ decreases the charge by $1$. Thus $(C(\mathfrak{g},x,f,k),d_{(0)})$ becomes a $\mathbb{Z}$-graded homology complex. Define \emph{the affine W-algebra} $W_{k}(\mathfrak{g},x,f)$ to be: as vector superspace the homology of this complex $W_{k}(\mathfrak{g},x,f):=H(C(\mathfrak{g},x,f,k),d_{(0)})$  together with the vertex algebra structure induced from $C(\mathfrak{g},x,f,k)$. The vertex algebra $W_{k}(\mathfrak{g},x,f)$ is also called the \emph{quantum reduction} associated to the quadruple $(\mathfrak{g},x,f,k)$. Define the Virasoro field of $C(\mathfrak{g},x,f,k)$ by 

\begin{equation}
L=L^{\mathfrak{g}}+\partial x+L^{ch}+L^{ne},\nonumber
\end{equation}

where 

\begin{equation}
L^{\mathfrak{g}}=\tfrac{1}{2(k+h^{\vee})}\sum_{\alpha\in S}(-1)^{p(u_{\alpha})}:u_{\alpha}u^{\alpha}:,\nonumber
\end{equation}

is given by the Sugawara construction, where $\{u^{\alpha}\}_{\alpha\in S}$ is the dual basis to $\{u_{\alpha}\}_{\alpha\in S}$, i.e., $(u_{\alpha}|u^{\beta})=\delta^{\beta}_{\alpha}$. Here we are assuming that $k\neq -h^{\vee}$, where $h^{\vee}$ denotes the dual Coxeter number of $\mathfrak{g}$.

\begin{equation}
L^{ch}=-\sum_{\alpha\in S_{+}}m_{\alpha}:\varphi^{\alpha}\partial\varphi_{\alpha}:+\sum_{\alpha\in S_{+}}(1-m_{\alpha}):(\partial\varphi^{\alpha})\varphi_{\alpha}:,\nonumber
\end{equation}

\begin{equation}
L^{ne}=\tfrac{1}{2}\sum_{\alpha\in S_{1/2}}:(\partial\Phi^{\alpha})\Phi_{\alpha}:.\nonumber
\end{equation}

The central charge of $L$ is given by

\begin{eqnarray} \label{centralchargeformula}
c(\mathfrak{g},x,f,k)&=&\frac{k\,sdim\, \mathfrak{g}}{k+h^{\vee}}-12k(x|x)\\
&&-\sum_{\alpha\in S_{+}}(-1)^{p(u_{\alpha})}\left(12 m_{\alpha}^{2}-12 m_{\alpha}+2\right)-\frac{1}{2} sdim\, \mathfrak{g}_{1/2}.\nonumber
\end{eqnarray}

With respect to $L$ the fields $a\,(a\in \mathfrak{g}_{j}),$ $\varphi_{\alpha}, \varphi^{\alpha}$ $(\alpha\in S_{+})$ and $\Phi_{\alpha}(\alpha\in S_{1/2})$ are primary vectors except for $a\,(a\in \mathfrak{g}_{0})$ such that $(a|x)\neq 0$, and the conformal weights are as follows: $\Delta(a)=1-j\; (a\in \mathfrak{g}_{j})$, $\Delta(\varphi_{\alpha})=1-m_{\alpha}$, $\Delta(\varphi^{\alpha})=m_{\alpha}$ and $\Delta(\Phi_{\alpha})=\frac{1}{2}$. In \cite{KacRoanWakimoto03} is proved that $d_{(0)}L=0$, then the homology class of $L$ (which does not vanish) defines the Virasoro field of $W_{k}(\mathfrak{g},x,f)$, which is again denoted by $L$.

To construct other fields of $W_{k}(\mathfrak{g},x,f)$ define for each $v\in\mathfrak{g}_{j}$

\begin{equation}
J^{(v)}=v+\sum_{\alpha,\beta\in S_{+}}(-1)^{p(u_{\alpha})}c^{\alpha}_{\beta}(v):\varphi_{\alpha}\varphi^{\beta}:,\nonumber
\end{equation}

where the numbers $c^{\alpha}_{\beta}(v)$ are given by $[v,u_{\alpha}]=\sum_{\alpha\in S}c^{\alpha}_{\beta}(v)u_{\alpha}$. The field $J^{(v)}\in C(\mathfrak{g},x,f,k) $ has the same charge, the same parity and the same conformal weight as the field $v$. The $\lambda$-bracket between these fields is as follows:

\begin{equation}\label{currentbracket}
[{J^{(v)}}_{\lambda}J^{(v')}]=J^{([v,v'])}+\lambda\left(k(v|v')+\tfrac{1}{2}\left(\kappa_{\mathfrak{g}}(v,v')-\kappa_{\mathfrak{g}_{0}}(v,v')\right)\right),
\end{equation}

if $v\in\mathfrak{g}_{i}, v'\in \mathfrak{g}_{j}$ and $ij\ge 0$ where $\kappa_{\mathfrak{g}}$(resp.$\kappa_{\mathfrak{g}_{0}}$) denotes the Killing form on $\mathfrak{g}$ (resp. $\mathfrak{g}_{0}$).

Denote by $C^{-}$ the vertex subalgebra of the vertex algebra $C(\mathfrak{g},x,f,k)$ generated by
the fields $J^{(u)}$ for all $u \in \mathfrak{g}_{\le}$, the fields $\varphi^{\alpha}$ for all $\alpha\in S_{+}$ and the fields $\Phi_{\alpha}$ for all $\alpha \in S_{1/2}$. One of the main theorems on the structure of the vertex algebra $W_{k}(\mathfrak{g},x,f)$ is the following:

\begin{theorem}\cite[Theorem 4.1]{KacWakimoto04}\label{maintheorem}
Let $\mathfrak{g}$ be a simple finite-dimensional Lie superalgebra with an invariant bilinear
form $(.|.)$ and let $x, f$ be a pair of even elements of $\mathfrak{g}$ such that $ad\;x$ is diagonalizable with eigenvalues in $\frac{1}{2}\mathbb{Z}$ and $[x,f]=-f$. Suppose that all eigenvalues of $ad\;x$ on $g^{f}$ (the centralizer of $f$) are non-positive: $g^{f}=\oplus_{j\le 0}g^{f}_{j}$. Then
\begin{itemize}
\item[a)] For each $a\in \mathfrak{g}^{f}_{-j}(j\ge 0)$ there exists a $d_{(0)}$-closed field $J^{\{a\}}$ in $C^{-}$ of conformal weight $1 + j$ (with respect to $L$) such that $J^{\{a\}}-J^{(a)}$ is a linear combination of normal ordered products of the fields $J^{(b)}$, where $b\in\mathfrak{g}_{-s}$, $0\le s< j$, the fields $\Phi_{\alpha}$, where $\alpha\in S_{1/2}$, and the derivatives of these fields.
\item[b)] The homology classes of the fields $J^{\{a_{i}\}}$, where $a_1, a_2,\dots$ is a basis of $g^{f}$ compatible with its $\frac{1}{2}\mathbb{Z}$-gradation, strongly generate the vertex algebra $W_{k}(\mathfrak{g}, x, f)$.
\item[c)]$H_{0}\left(C(\mathfrak{g},x,f,k),d_{(0)}\right)=W_{k}(\mathfrak{g},x,f)$ and $H_{j}\left(C(\mathfrak{g},x,f,k),d_{(0)}\right)=0$  if $j\neq 0$.
\end{itemize}

\end{theorem}

\begin{remark}
The complex $\left(C(\mathfrak{g},x,f,k),d_{(0)}\right)$ is formal, that is, the vertex algebra $W_{k}(\mathfrak{g},x,f)$ is a subalgebra of $C(\mathfrak{g},x,f,k)$ consisting of $d_{(0)}$-closed charge $0$ elements of $C^{-}$, furthermore the $J^{\{a\}}$ can be computed recursively, for example in the case $a\in g^{f}_{-1/2}$ the solution is unique an is given by:
\end{remark}

\begin{theorem}\cite[Theorem 2.1 (d)]{KacWakimoto04}\label{reconstruction of the fields}

For $v \in \mathfrak{g}_{-1/2}$ let 

\begin{eqnarray}
G^{\{v\}}&=& J^{(v)}+\sum_{\beta\in S_{1/2}}:J^{([v,u_{\beta}])}\Phi^{\beta}:+\tfrac{(-1)^{p(v)+1}}{3}\sum_{\alpha,\beta\in S_{1/2}}:\Phi^{\alpha}\Phi^{\beta}\Phi_{[u_{\beta}[u_{\alpha},v]]}:\nonumber\\
&&-\sum_{\beta\in S_{1/2}}\left(k(v|u_{\beta})+str_{g_{+}}(ad\,v)(ad\,u_{\beta})\right)\partial\Phi^{\beta}, \nonumber
\end{eqnarray}
Then provided that  $v \in \mathfrak{g}^{f}_{-1/2}$, we have $d_{(0)}(G^{\{v\}})=0$, hence the homology class of $G^{\{v\}}$ defines a field of the vertex algebra $W_{k}(\mathfrak{g},x,f)$ of conformal weight $\frac{3}{2}$. This field is primary.
\end{theorem}

\begin{remark}\label{W algebra as a subalegbra}
In the case $\mathfrak{g}^{f}\subset g_{\le}$ Theorem \ref{maintheorem} and the identity (\ref{currentbracket}) provides a construction of the vertex algebra $W_{k}(\mathfrak{g},x,f)$ as a subalgebra of $V_{\nu_{k}}(\mathfrak{g}_{\le})\otimes F(\mathfrak{g}_{1/2})$ where $\nu_{k}$ is the 2-cocycle on $\mathfrak{g}_{\le}[t,t^{-1}]$ given by 
\begin{equation}\label{newcocycle}
\nu_{k}(at^{m},bt^{n})=m\delta_{m,-n}\left(k(a|b)+\tfrac{1}{2}\left(\kappa_{\mathfrak{g}}(a,b)-\kappa_{\mathfrak{g}_{0}}(a,b)\right)\right)
\end{equation}
for $a,b\in \mathfrak{g_{\le}}$ and $m,n\in \mathbb{Z}$.
\end{remark}

\begin{remark}\label{free field realization of W}
Furthermore if this 2-cocycle is trivial outside $\mathfrak{g}_{0}[t,t^{-1}]$, the cannonical homomorphism $\mathfrak{g}_{\le}\rightarrow\mathfrak{g}_{0}$ induces a homomorphism from $V_{\nu_{k}}(\mathfrak{g}_{\le})\otimes F(\mathfrak{g}_{1/2})$ to $V_{\nu_{k}}(\mathfrak{g}_{0})\otimes F(\mathfrak{g}_{1/2})$, obtaining in this way a free field realization of $W_{k}(\mathfrak{g},x,f)$  inside $V_{\nu_{k}}(\mathfrak{g}_{0})\otimes F(\mathfrak{g}_{1/2})$.
\end{remark}

\section{Quantum Hamiltonian Reduction of $D(2,1;\alpha)$}\label{quantum reduction of the algebra}

In this section we prove that the family $SW(\frac{3}{2},\frac{3}{2},2)$ of W-algebras which has generators $\{G, H, L, \tilde{M}, W, U\}$ of conformal weights $\left(\frac{3}{2}, \frac{3}{2}, 2, 2, 2,\frac{5}{2}\right)$ and relations as given in Appendix \ref{App:AppendixB} can be obtained as the quantum Hamiltonian reduction of $D(2,1;\alpha)$. As a collorary we obtain a free field realization of this family. As a particular case we obtain a free-field realization of the Shatashvili-Vafa $G_{2}$ algebra on a space of three free Bosons and three free Fermions.

The Lie superalgebra $D(2,1;\alpha)$ where $\alpha\in\mathbb{C}\backslash \{-1,0\}$ is a one-parameter family of exceptional Lie superalgebras of rank $3$ and dimension $17$, which contains $D(2,1)=osp(4,2)$ as special cases (when $\alpha\in\{1,-\frac{1}{2},-2\}$), see \cite{Kac77}.

We present $\mathfrak{g}=D(2,1;\alpha)$ as the contragradient Lie superalgebra associated to the Cartan matrix $A=(a_{ij})_{i,j}$ and $\tau=\{1,2,3\}$

\begin{eqnarray}\label{Cartanmatrix}
{(a_{ij})}_{i,j=1}^{3}= \left( \begin{array}{ccc}
0 & 1 & \alpha \\
1 & 0 & -1-\alpha \\
\alpha & -1-\alpha & 0 \end{array} \right).
\end{eqnarray}

We have generators $\{h_{1},h_{2},h_{3},e_{1},e_{2},e_{3},f_{1},f_{2},f_{3}\}$, $h_i$ being even for all $i$ and $e_i,f_i$ being odd for all $i$ and relations

\begin{equation}
[e_i,f_j]=\delta_{ij}h_{i},\;\; [h_i,e_j]=a_{ij}e_j,\;\; [h_i,f_j]=-a_{ij}f_j.\nonumber
\end{equation}

Introduce the elements:

\begin{equation}
[e_1,e_2]=:e_{12},\;\; [e_1,e_3]=:e_{13},\;\; [e_2,e_3]=:e_{23},\;\; [e_1,e_{23}]=:e_{123},\nonumber
\end{equation}

\begin{equation}
[f_1,f_2]=:f_{12},\;\; [f_1,f_3]=:f_{13},\;\; [f_2,f_3]=:f_{23},\;\; [f_1,f_{23}]=:f_{123}.\nonumber
\end{equation}

Recall that $\mathfrak{g}$ has vanishing Killing form and consequently the dual Coxeter number $h^{\vee}=0$. Fix the following non-degenerate even supersymmetric invariant bilinear form $(.|.)$

\begin{equation}
(h_i,h_j)=a_{ij},\;\; (e_i,f_j)=\delta_{ij},\;\; (e_{12},f_{12})=(f_{12},e_{12})=-1,\nonumber
\end{equation}
\begin{equation}
(e_{13},f_{13})=(f_{13},e_{13})=-\alpha,\;\; (e_{23},f_{23})=(f_{23},e_{23})=1+\alpha,\nonumber
\end{equation}
\begin{equation}
(e_{123},f_{123})=-(f_{123},e_{123})=(1+\alpha)^2.\nonumber
\end{equation}

To perform the quantum Hamiltonian reduction we take the pair $(x,f)$:
\begin{eqnarray}
x:=\tfrac{(\alpha+1)}{2\alpha}h_1+\tfrac{\alpha}{2(\alpha+1)}h_2+\tfrac{1}{2\alpha(\alpha+1)}h_3,&f:=f_{12}+f_{13}+f_{23}.\nonumber
\end{eqnarray}
This pair together with $e=(-\tfrac{1}{2})e_{12}+(-\tfrac{1}{2\alpha^{2}})e_{13}+(-\tfrac{1}{2(\alpha+1)^{2}})e_{23}$ forms an $sl_{2}$ triple. We have the following eigenspace decomposition of the algebra with respect to $ad$ $x$:
\[  \begin{array}{ccccccc}
\mathfrak{g}_{-3/2} & \mathfrak{g}_{-1} & \mathfrak{g}_{-1/2} & \mathfrak{g}_{0} & \mathfrak{g}_{1/2} & \mathfrak{g}_{1} & \mathfrak{g}_{3/2} \\
f_{123} & f_{12} & f_{1} & h_{1} & e_{1} & e_{12} & e_{123} \\
 & f_{13} & f_{2} & h_{2} & e_{2} & e_{13} & \\
 & f_{23} & f_{3} & h_{3} & e_{3} & e_{23} & \\ 
 \end{array} \]
 
Furthermore $\mathfrak{g}^{f}=\mathfrak{g}_{-1/2}^{f}\oplus \mathfrak{g}_{-1}^{f}\oplus \mathfrak{g}_{-3/2}^{f}$ with dim $ \mathfrak{g}_{-1/2}^{f}=2$ , dim $ \mathfrak{g}_{-1}^{f}=3$ and dim $\mathfrak{g}_{-3/2}^{f}=1$. This shows that the algebra $W_{k}(\mathfrak{g},x,f)$ has six generators with the expected conformal weights.

The set of vectors  $\{e_{1},e_{2},e_{3}\}$ is a basis of $\mathfrak{g}_{1/2}$ , denote by $\Phi_{1}:=e_{1}$, $\Phi_{2}:=e_{2}$ and $\Phi_{3}:=e_{3}$ the corresponding free neutral Fermions. The non-zero values of the (symmetric) bilinear form $\left<.|.\right>_{ne}$ on $\mathfrak{g}_{1/2}$ are given by:

$$\left<\Phi_{1}|\Phi_{2}\right>_{ne}=-1,\;\;\left<\Phi_{1}|\Phi_{3}\right>_{ne}=-\alpha,\;\; \left<\Phi_{2}|\Phi_{3}\right>_{ne}=1+\alpha,$$

\noindent
(note that this is exactly minus the Cartan matrix of $D(2,1;\alpha)$.) Then the free neutral fermions satisfy the following non-zero $\lambda$-brackets:

$$[{\Phi_{1}}_{\lambda}\Phi_{2}]=-1,\;\;[{\Phi_{1}}_{\lambda}\Phi_{3}]=-\alpha,\;\;[{\Phi_{2}}_{\lambda}\Phi_{3}]=1+\alpha,$$

\noindent
and the dual free neutral fermions with respect to $\left<.|.\right>_{ne}$ are:

\begin{eqnarray}
\Phi^{1}&=&(-\tfrac{1+\alpha}{2\alpha})\Phi_{1}+(-\tfrac{1}{2})\Phi_{2}+(-\tfrac{1}{2\alpha})\Phi_{3}, \nonumber\\ \Phi^{2}&=&(-\tfrac{1}{2})\Phi_{1}+(-\tfrac{\alpha}{2+2\alpha})\Phi_{2}+(\tfrac{1}{2+2\alpha})\Phi_{3}, \nonumber\\ \Phi^{3}&=&(-\tfrac{1}{2\alpha})\Phi_{1}+(\tfrac{1}{2+2\alpha})\Phi_{2}+(-\tfrac{1}{2\alpha+2\alpha^{2}})\Phi_{3}.\nonumber
\end{eqnarray}

\noindent
We fix the basis $\{h_{1},h_{2},h_{3},f_{1},f_{2},f_{3},f_{12},f_{13},f_{23},f_{123}\}$  of $\mathfrak{g}_{\le}$ compatible with the $\frac{1}{2}\mathbb{Z}$ and $\mathbb{Z}_{2}$ gradation of $\mathfrak{g}$. We consider the building blocks $J^{(v)}$ for each $v$ that belongs to the above basis, (\ref{currentbracket}) reduces to
\begin{equation}
[{J^{(v)}}_{\lambda}J^{(v')}]=J^{([v,v'])}+\lambda k(v|v'),\nonumber
\end{equation}
\noindent
because the Killing form  $\kappa_{\mathfrak{g}}$  of $\mathfrak{g}$ is zero and $\mathfrak{g}_{0}$ equals the Cartan subalgebra $\mathfrak{h}$ of $\mathfrak{g}$, that is, the generators $J^{(v)}$ obey the same commutation relations as the generators of $V_{k}(D(2,1;\alpha))$. Using Remark \ref{W algebra as a subalegbra} we obtain that $W_{k}(\mathfrak{g},x,f)$ is a subalgebra of $V_{k}(\mathfrak{g}_{\le})\otimes F(\mathfrak{g}_{1/2})$. For this reason and to simplify the notation we denote $J^{(v)}$ simply by $v$. Furthermore as the cocycle (\ref{newcocycle}) is the original cocycle of $V_{k}(D(2,1;\alpha))$ and this cocycle is trivial in $\mathfrak{g}_{\le}$ outside $\mathfrak{g}_{0}=\mathfrak{h}$, Remark \ref{free field realization of W} gives a free field realization of $W_{k}(\mathfrak{g},x,f)$ inside $V_{k}(\mathfrak{h})\otimes F(\mathfrak{g}_{1/2})$.

 Let $J^{\{f_i\}}$ denote the $d_{0}$-closed fields associated to $\{f_{i}\}_{i=1}^{3}$ provided by Theorem \ref{maintheorem}. Using Theorem \ref{reconstruction of the fields} we can compute $J^{\{f_i\}}$ explicitly: 
  
\begin{eqnarray}
J^{\{f_1\}}&=&f_{1}+\left(\tfrac{\alpha^{2}-1}{3}\right):\Phi^{1}\Phi^{2}\Phi^{3}:+:\Phi^{1}h_{1}:+k \partial\Phi^{1},\nonumber \\
J^{\{f_2\}}&=&f_{2}+\left(\tfrac{\alpha(\alpha+2) }{3}\right):\Phi^{1}\Phi^{2}\Phi^{3}:+:\Phi^{2}h_{2}:+k \partial\Phi^{2},\nonumber \\
J^{\{f_3\}}&=&f_{3}+\left(\tfrac{2\alpha+1}{3}\right):\Phi^{1}\Phi^{2}\Phi^{3}:+:\Phi^{3}h_{3}:+k \partial\Phi^{3}. \nonumber
\end{eqnarray}

We can compute the other fields $J^{\{f_{1,2}\}},J^{\{f_{1,3}\}},J^{\{f_{2,3}\}},J^{\{f_{1,2,3}\}}$ given by Theorem \ref{maintheorem} that jointly with $\{J^{\{f_i\}}\}_{i=1}^{3}$ strongly generate $W_{k}(\mathfrak{g},x,f)$, but in the $SW(\frac{3}{2},\frac{3}{2},2)$ superconformal algebra we can recover (using $\lambda$-brackets) all the fields from the generators in conformal weight $\frac{3}{2}$, i.e., $G$ and $H$ (see Appendix \ref{App:AppendixB}). Thus we only need to construct $G$ and $H$ from $\{J^{\{f_i\}}\}_{i=1}^{3}$.

In order to do that observe that:
\begin{equation}\label{condition to be conformal weight 2}
a_{1}f_{1}+a_{2}f_{2}+a_{3}f_{3}\in \mathfrak{g}^{f}_{-1/2} \Leftrightarrow a_{1}+a_{2}(-\tfrac{\alpha}{\alpha+1})+a_{3}(-\tfrac{1}{\alpha+1})=0,
\end{equation}

\noindent
and that the central charge of the Virasoro field of $W_{k}(\mathfrak{g},x,f)$ given by formula (\ref{centralchargeformula}) is $c(\alpha,k)=\tfrac{9}{2}-12k(x|x)=\tfrac{9}{2}-\tfrac{6 k \left(1+\alpha +\alpha ^2\right)}{\alpha(1+\alpha )}$.

We want to define a field $G$ such that $\{G,L:=\tfrac{1}{2}{G}_{(0)}G\}$ generate an $N=1$ superconformal algebra with the above central charge, this is accomplished taking $a_{1}=a_{2}=a_{3}=\tfrac{i}{\sqrt{k}}$, i.e.,

\begin{equation}
G:=\tfrac{i}{\sqrt{k}}\left(J^{\{f_1\}}+J^{\{f_2\}}+J^{\{f_3\}}\right).\nonumber
\end{equation}

We are looking for a vector $H$ of conformal weight $\frac{3}{2}$, such that:

\begin{equation}\label{definition of H hat}
{G}_{(j)}H=0, \;\;j>0,
\end{equation}

\noindent
The most general vector of conformal weight $\frac{3}{2}$ given by (\ref{condition to be conformal weight 2}) is 

\begin{equation}
\left(\tfrac{\alpha}{\alpha+1}a_{2}+\tfrac{1}{\alpha+1}a_{3}\right)J^{\{f_1\}}+a_{2}J^{\{f_2\}}+a_{3} J^{\{f_3\}},\nonumber
\end{equation}

\noindent
(\ref{definition of H hat}) imposes the condition $a_{2}\alpha(-1+2 k)(1+2 \alpha )+ a_{3}(2 k-\alpha )(2+\alpha )=0,$ which has as solution
\begin{eqnarray}
{a_{1}}'&:=&\alpha(-1+\alpha )(1+2k+\alpha ),\nonumber\\
{a_{2}}'&:=&(-1)(2k-\alpha )(2+\alpha )(1+\alpha ),\nonumber\\
{a_{3}}'&:=&\alpha(-1+2k)(1+2\alpha)(1+\alpha).\nonumber
\end{eqnarray}
It follows from ${H}_{(2)}H=\frac{2c}{3}$ (cf. (\ref{bracketHwithH})) that we need to rescale this solution to define $H=\sum_{i=1}^{3}a_{i}J^{\{f_i\}}$ with
\begin{equation}
a_{i}:=\left(-\tfrac{3}{2} (-1+2 k)\alpha ^2 (1+\alpha )^2 \left(2 k+4 k^2-\alpha  (1+\alpha )\right)\right)^{-\tfrac{1}{2}}{a_{i}}'.\nonumber
\end{equation}

We can obtain all other generators from $G$ and $H$, to perform this computations we use Thielemans's software \cite{Thielemans91}. Listed below are the explicit expressions of all the generators of $W_{k}(\mathfrak{g},x,f)$ as a subalgebra of $V_{k}(\mathfrak{g}_{\le})\otimes F(\mathfrak{g}_{1/2})$:

\begin{eqnarray}
G&=&\tfrac{i}{\sqrt{k}}f_{1}+\tfrac{i}{\sqrt{k}}f_{2}+\tfrac{i}{\sqrt{k}}f_{3}+\tfrac{i}{\sqrt{k}}:\Phi^{1}h_{1}:+\tfrac{i}{\sqrt{k}}:\Phi^{2}h_{2}:+\tfrac{i}{\sqrt{k}}:\Phi^{3}h_{3}:\nonumber\\
&&+i\sqrt{k} \partial\Phi^{1}+i\sqrt{k} \partial\Phi^{2}+i\sqrt{k} \partial\Phi^{3},\nonumber
\end{eqnarray}

\begin{eqnarray}
L&=&-\tfrac{1}{k}f_{12}-\tfrac{1}{k}f_{13}-\tfrac{1}{k}f_{23}+\tfrac{(1+\alpha )}{4 k \alpha }:h_{1}h_{1}:+\tfrac{1}{2k}:h_{1}h_{2}:+\tfrac{1}{2k\alpha}:h_{1}h_{3}:\nonumber\\
&&+\tfrac{\alpha}{4 k+4 k \alpha }:h_{2}h_{2}:-\tfrac{1}{2 k+2 k \alpha }:h_{2}h_{3}:+\tfrac{1}{4 k \alpha +4 k \alpha ^2}:h_{3}h_{3}:+\tfrac{1}{k}:\Phi^{1}f_{2}:\nonumber\\
&&+\tfrac{\alpha}{k}:\Phi^{1}f_{3}:+\tfrac{1}{2}:\Phi^{1}\partial\Phi^{2}:+\tfrac{1}{2} \alpha:\Phi^{1}\partial\Phi^{3}:+\tfrac{1}{k}:\Phi^{2}f_{1}:-\tfrac{(1+\alpha )}{k}:\Phi^{2}f_{3}:\nonumber\\
&&+\tfrac{1}{2} (-1-\alpha ):\Phi^{2}\partial\Phi^{3}:+\tfrac{\alpha}{k}:\Phi^{3}f_{1}:-\tfrac{(1+\alpha )}{k} :\Phi^{3}f_{2}:-\tfrac{1}{2}:\partial\Phi^{1}\Phi^{2}:\nonumber\\
&&-\tfrac{1}{2} \alpha :\partial\Phi^{1}\Phi^{3}:+\tfrac{1}{2} (1+\alpha ):\partial\Phi^{2}\Phi^{3}:+\tfrac{(1+\alpha)}{2\alpha }\partial h_{1}+\tfrac{\alpha}{2+2 \alpha }\partial h_{2}\nonumber\\
&&+\tfrac{1}{2\alpha+2\alpha^2}\partial h_{3},\nonumber
\end{eqnarray}

\small
\begin{eqnarray}
H&=&\tfrac{1}{\sqrt{-\tfrac{3}{2} (-1+2 k)\alpha ^2 (1+\alpha )^2 \left(2 k+4 k^2-\alpha  (1+\alpha )\right)}}\left(\left(-1+\alpha\right) \alpha  \left(1+2k+\alpha\right) \left(f_{1}+ :\Phi^{1}h_{1}:\right.\right.\nonumber\\
&&\left.\left.+k \partial\Phi^{1}\right)-(2k-\alpha ) \left(2+3 \alpha +\alpha ^2\right) (f_{2}+ :\Phi^{2}h_{2}:+k\partial\Phi^{2})\right.\nonumber\\
&&\left.+(-1+2 k) \alpha  \left(1+3 \alpha +2 \alpha ^2\right) (f_{3}+ :\Phi^{3}h_{3}:+k\partial\Phi ^{3})\right.\nonumber\\
&&\left.+\alpha  (1+\alpha ) \left(-3 \alpha  (1+\alpha )+4 k \left(1+\alpha +\alpha ^2\right)\right) :\Phi^{1}\Phi^{2}\Phi^{3}:\right),\nonumber
\end{eqnarray}
\normalsize

\small
\begin{eqnarray}
\tilde{M}&=&\tfrac{1}{\sqrt{-\tfrac{3}{2} (-1+2 k)\alpha ^2 (1+\alpha )^2 \left(2 k+4 k^2-\alpha  (1+\alpha )\right)}}\left(\tfrac{i (1+2 \alpha ) \left(-4 k+\alpha +\alpha ^2\right) }{\sqrt{k}}\left(f_{12}-\tfrac{1}{2}:h_{1}h_{2}:\right.\right.\nonumber\\
&&\left.\left.- :\Phi^{1}f_{2}- :\Phi^{2}f_{1}:\right)+\tfrac{i (2+\alpha ) (-1+(-1+4 k) \alpha ) }{\sqrt{k}}\left(\alpha  f_{13}-\tfrac{1}{2}:h_{1}h_{3}:-\alpha ^2:\Phi^{1}f_{3}:\right.\right.\nonumber\\
&&\left.\left.- \alpha^2:\Phi^{3}f_{1}:\right)+\tfrac{i \left(-1+\alpha ^2\right) (-\alpha +4 k (1+\alpha )) }{\sqrt{k}}\left(f_{23}+\tfrac{1}{2 (1+\alpha )}:h_{2}h_{3}:\right.\right.\nonumber\\
&&\left.\left.+ (1+\alpha ):\Phi^{3}f_{2}:+ (1+\alpha ):\Phi^{2}f_{3}:\right)\right.\nonumber\\
&&\left.-i \sqrt{k} (-1+\alpha ) (1+\alpha ) (1+2 k+\alpha )\left( \partial h_{1}+\tfrac{1}{2k}:h_{1}h_{1}:\right)\right.\nonumber\\
&&\left.+i \sqrt{k} (2 k-\alpha ) \alpha  (2+\alpha )\left( \partial h_{2}+\tfrac{1}{2k}:h_{2}h_{2}:\right)\right.\nonumber\\
&&\left.-i \sqrt{k} (-1+2 k) (1+2 \alpha )\left( \partial h_{3}+\tfrac{1}{2k}:h_{3}h_{3}:\right)\right.\nonumber\\
&&\left.+\tfrac{i\left(-3 \alpha  (1+\alpha )+4 k \left(1+\alpha +\alpha ^2\right)\right) }{2 \sqrt{k}}\left(-(1+\alpha ):\Phi^{1}\Phi^{2}h_{1}:+\alpha:\Phi^{1}\Phi^{2}h_{2}:-:\Phi^{1}\Phi^{2}h_{3}:\right.\right.\nonumber\\
&&\left.\left.+\alpha(1+\alpha ):\Phi^{1}\Phi^{3}h_{1}:+\alpha^2:\Phi^{1}\Phi^{3}h_{2}:-\alpha:\Phi^{1}\Phi^{3}h_{3}:
-(1+\alpha )^2:\Phi^{2}\Phi^{3}h_{1}:\right.\right.\nonumber\\
&&\left.\left.-\alpha(1+\alpha ):\Phi^{2}\Phi^{3}h_{2}:-(1+\alpha ):\Phi^{2}\Phi^{3}h_{3}:\right)\right.\nonumber\\
&&\left.-i \sqrt{k} (-1+\alpha ) \alpha  (1+2 k+\alpha ):\Phi^{1}\partial\Phi^{2}:\right.\nonumber\\
&&\left.-i \sqrt{k} (-1+\alpha ) \alpha ^2 (1+2 k+\alpha ):\Phi^{1}\partial\Phi^{3}:\right.\nonumber\\
&&\left.-i \sqrt{k} (2 k-\alpha ) (1+\alpha )^2 (2+\alpha ) :\Phi^{2}\partial\Phi^{3}:\right.\nonumber\\
&&\left.-i \sqrt{k} (2 k-\alpha ) \left(2+3 \alpha +\alpha ^2\right) :\partial\Phi^{1}\Phi^{2}:\right.\nonumber\\
&&\left.+i \sqrt{k} (-1+2 k) \alpha ^2 \left(1+3 \alpha +2 \alpha ^2\right) :\partial\Phi^{1}\Phi^{3}:\right.\nonumber\\
&&\left.-i \sqrt{k} (-1+2 k) \alpha  (1+\alpha )^2 (1+2 \alpha ) :\partial\Phi^{2}\Phi^{3}:\right)\nonumber,
\end{eqnarray}
\normalsize

\small
\begin{eqnarray}
W&=&\tfrac{\mu }{\left(-3 \alpha  (1+\alpha )+4 k \left(1+\alpha +\alpha ^2\right)\right)}\left(\tfrac{\left(2 k+\alpha +\alpha ^2\right)}{k}\left(- f_{12}+\tfrac{1}{2 } :h_{1}h_{2}:+ :\Phi^{1}f_{2}:\right.\right.\nonumber\\
&&\left.\left.+ :\Phi^{2}f_{1}:\right)+\tfrac{(1+\alpha +2 k \alpha ) }{k}\left(-\alpha  f_{13}+\tfrac{1}{2 } :h_{1}h_{3}:+\alpha ^2 :\Phi^{1}f_{3}:+\alpha ^2:\Phi^{3}f_{1}:\right)\right.\nonumber\\
&&\left.+\tfrac{(1+\alpha ) (\alpha +2 k (1+\alpha )) }{k}\left(-f_{23}-(1+\alpha ):\Phi^{2}f_{3}:-(1+\alpha ):\Phi^{3}f_{2}\right)\right.\nonumber\\
&&\left.+\tfrac{(1+\alpha ) (1+2 k+\alpha )}{4 k}:h_{1}h_{1}:+\tfrac{\alpha  (-2 k+\alpha )}{4k}:h_{2}h_{2}:+\tfrac{-2k+\alpha(-2k-1)}{2k}:h_{2}h_{3}:\right.\nonumber\\
&&\left.+\tfrac{-2k+1}{4k}:h_{3}h_{3}:+\left(-1+\alpha ^2\right):\Phi^{1}\Phi^{2}h_{1}:-\alpha  (2+\alpha ):\Phi^{1}\Phi^{2}h_{2}:\right.\nonumber\\
&&\left.+(-1-2 \alpha ) :\Phi^{1}\Phi^{2}h_{3}:+\left(\alpha -\alpha ^3\right) :\Phi^{1}\Phi^{3}h_{1}:-\alpha ^2 (2+\alpha ):\Phi^{1}\Phi^{3}h_{2}:\right.\nonumber\\
&&\left.-\alpha  (1+2 \alpha ):\Phi^{1}\Phi^{3}h_{3}:-\alpha  (1+2 k+\alpha ) :\Phi^{1}\partial\Phi^{2}:\right.\nonumber\\
&&\left.-\alpha ^2 (1+2 k+\alpha ):\Phi^{1}\partial\Phi^{3}:+(-1+\alpha ) (1+\alpha )^2 :\Phi^{2}\Phi^{3}h_{1}:
\right.\nonumber\\
&&\left.+\alpha  \left(2+3 \alpha +\alpha ^2\right) :\Phi^{2}\Phi^{3}h_{2}:+\left(-1-3 \alpha -2 \alpha ^2\right):\Phi^{2}\Phi^{3}h_{3}:\right.\nonumber\\
&&\left.-(2 k-\alpha ) (1+\alpha )^2 :\Phi^{2}\partial\Phi^{3}:-(2 k-\alpha ) (1+\alpha ):\partial\Phi^{1}\Phi^{2}:\right.\nonumber\\
&&\left.-(-1+2 k) \alpha ^2 (1+\alpha ) :\partial\Phi^{1}\Phi^{3}:+(-1+2 k) \alpha  (1+\alpha )^2 :\partial\Phi^{2}\Phi^{3}:\right.\nonumber\\
&&\left.+\tfrac{1}{2} (1+\alpha ) (1+2 k+\alpha ) \partial h_{1}+\tfrac{1}{2} \alpha  (-2 k+\alpha ) \partial h_{2}+\tfrac{1-2k}{2}\partial h_{3}\right),\nonumber
\end{eqnarray}
\normalsize

\small
\begin{align*}
U&=\tfrac{\mu }{\left(-3 \alpha  (1+\alpha )+4 k \left(1+\alpha +\alpha ^2\right)\right)}\left(-\tfrac{6i\alpha}{\sqrt{k}}f_{123}+\tfrac{3i(1+\alpha)}{\sqrt{k}}:h_{1}f_{2}:+\tfrac{3i\alpha(1+\alpha ) }{\sqrt{k}}:h_{1}f_{3}:\right.\\
&\quad\left.-\tfrac{3 i \alpha}{\sqrt{k}}:h_{2}f_{1}:+\tfrac{3 i \alpha  (1+\alpha ) }{\sqrt{k}}:h_{2}f_{3}:-\tfrac{3 i \alpha}{\sqrt{k}}:h_{3}f_{1}:\right.\\
&\quad\left.+\frac{3 i (1+\alpha ) }{\sqrt{k}}:h_{3}f_{2}:+\tfrac{6 i \alpha  (1+\alpha )}{\sqrt{k}}:\Phi^{1}f_{23}:-\tfrac{3 i \alpha  }{\sqrt{k}}:\Phi^{1}h_{1}h_{2}:\right.\\
&\quad\left.-\tfrac{3 i \alpha }{\sqrt{k}}:\Phi^{1}h_{1}h_{3}:-\tfrac{3 i \alpha  (1+\alpha )}{\sqrt{k}}:\Phi^{1}\Phi^{2}f_{1}:+\tfrac{3 i \alpha  (1+\alpha ) }{\sqrt{k}}:\Phi^{1}\Phi^{2}f_{2}:\right.\\
&\quad\left.+\tfrac{3 i \alpha  \left(1+3 \alpha +2 \alpha ^2\right) }{\sqrt{k}}:\Phi^{1}\Phi^{2}f_{3}:+i \sqrt{k} \alpha  \left(1+3 \alpha +2 \alpha ^2\right) :\Phi^{1}\Phi^{2}\partial\Phi^{3}:\right.\\
&\quad\left.-\tfrac{3 i \alpha ^2 (1+\alpha )}{\sqrt{k}}:\Phi^{1}\Phi^{3}f_{1}:+\tfrac{3 i \alpha  \left(2+3 \alpha +\alpha ^2\right) }{\sqrt{k}}:\Phi^{1}\Phi^{3}f_{2}:+\tfrac{3 i \alpha ^2 (1+\alpha )}{\sqrt{k}}:\Phi^{1}\Phi^{3}f_{3}:\right.\\
&\quad\left.-3 i \sqrt{k} \alpha  (1+\alpha ):\Phi^{1}\partial\Phi^{2}\Phi^{2}:-i \sqrt{k} \alpha  \left(2+3 \alpha +\alpha ^2\right):\Phi^{1}\partial\Phi^{2}\Phi^{3}:\right.\\
&\quad\left.-3 i \sqrt{k} \alpha ^2 (1+\alpha ):\Phi^{1}\partial\Phi^{3}\Phi^{3}:-\tfrac{i \alpha  (1+2 k+\alpha )}{\sqrt{k}}:\Phi^{1}\partial h_{1}:\right.\\
&\quad\left.+\tfrac{6 i \alpha  (1+\alpha )}{\sqrt{k}}:\Phi^{2}f_{13}:+\tfrac{3 i (1+\alpha ) }{\sqrt{k}}:\Phi^{2}h_{1}h_{2}:+\tfrac{3 i (1+\alpha )}{\sqrt{k}}:\Phi^{2}h_{2}h_{3}:\right.\\
&\quad\left.-\tfrac{3 i \alpha  \left(-1+\alpha ^2\right) }{\sqrt{k}}:\Phi^{2}\Phi^{3}f_{1}:+\tfrac{3 i \alpha  (1+\alpha )^2 }{\sqrt{k}}:\Phi^{2}\Phi^{3}f_{2}:-\tfrac{3 i \alpha  (1+\alpha )^2 }{\sqrt{k}}:\Phi^{2}\Phi^{3}f_{3}:\right.\\
&\quad\left.+3 i \sqrt{k} \alpha  (1+\alpha )^2 :\Phi^{2}\partial\Phi^{3}\Phi^{3}:+\tfrac{i (2 k-\alpha ) (1+\alpha ) }{\sqrt{k}}:\Phi^{2}\partial h_{2}:+\tfrac{6 i \alpha  (1+\alpha )}{\sqrt{k}}:\Phi^{3}f_{12}:\right.\\
&\quad\left.+\tfrac{3 i \alpha  (1+\alpha )}{\sqrt{k}}:\Phi^{3}h_{1}h_{3}:+\tfrac{3 i \alpha  (1+\alpha )}{\sqrt{k}}:\Phi^{3}h_{2}h_{3}:+\tfrac{i (-1+2 k) \alpha  (1+\alpha )}{\sqrt{k}}:\Phi^ {3}\partial h_{3}:\right.\\
&\quad\left.-\tfrac{2 i (-1+k-\alpha ) \alpha  }{\sqrt{k}}:\partial\Phi^{1}h_{1}:-3 i \sqrt{k} \alpha :\partial\Phi^{1}h_{2}:-3 i \sqrt{k} \alpha:\partial\Phi^{1}h_{3}:\right.\\
&\quad\left.-3 i \sqrt{k} \alpha  (1+\alpha ):\partial\Phi^{1}\Phi^{1}\Phi^{2}:-3 i \sqrt{k} \alpha ^2 (1+\alpha ):\partial\Phi^{1}\Phi^{1}\Phi^{3}:\right.\\
&\quad\left.-i \sqrt{k} \alpha  \left(-1+\alpha ^2\right):\partial\Phi^{1}\Phi^{2}\Phi^{3}:+3 i \sqrt{k} (1+\alpha ):\partial\Phi^{2}h_{1}:\right.\\
&\quad\left.+\tfrac{2 i (1+\alpha ) (k+\alpha )}{\sqrt{k}}:\partial\Phi^{2}h_{2}:+3 i \sqrt{k} (1+\alpha ):\partial\Phi^{2}h_{3}:\right.\\
\displaybreak
&\quad\left.+3 i \sqrt{k} \alpha  (1+\alpha )^2:\partial\Phi^{2}\Phi^{2}\Phi^{3}:+3 i \sqrt{k} \alpha  (1+\alpha ):\partial\Phi^{3}h_{1}:\right.\\
&\quad\left.+3 i \sqrt{k} \alpha  (1+\alpha ):\partial\Phi^{3}h_{2}:+\tfrac{2 i (1+k) \alpha  (1+\alpha )}{\sqrt{k}}:\partial\Phi^{3}h_{3}:\right.\\
&\quad\left.-\tfrac{i \alpha  (1+2 k+\alpha )}{\sqrt{k}}\partial f_{1}+\tfrac{i (2 k-\alpha ) (1+\alpha )}{\sqrt{k}}\partial f_{2}+\tfrac{i (-1+2 k) \alpha  (1+\alpha ) }{\sqrt{k}}\partial f_{3}\right.\\
&\quad\left.-\tfrac{1}{2} i \sqrt{k} (-1+4 k-\alpha ) \alpha  \partial^2\Phi^{1}+\frac{1}{2} i \sqrt{k} (1+\alpha ) (4 k+\alpha ) \partial^2\Phi^{2}\right.\\
&\quad\left.+\frac{1}{2} i \sqrt{k} (1+4 k) \alpha  (1+\alpha ) \partial^2\Phi^{3}\right),
\end{align*}
\normalsize
where $\mu=\sqrt{\tfrac{9c(4+\varepsilon^{2})}{2(27-2c)}}$ and $\varepsilon(\alpha,k)=-\tfrac{4 i \sqrt{\frac{2}{3}} k^{3/2} (1+2 \alpha ) \left(-2+\alpha +\alpha ^2\right)}{3\sqrt{-(-1+2 k) \alpha ^2 (1+\alpha )^2 \left(2 k+4 k^2-\alpha  (1+\alpha )\right)}}$.

One can check straightforwardly with the aid of \cite{Thielemans91} that the $\lambda$-brackets of the algebra $W_{k}(\mathfrak{g},x,f)$ coincides with the $\lambda$-brackets of the family of superconformal algebras  $SW(\frac{3}{2},\frac{3}{2},2)$ with parameters $\left(c(\alpha,k),\varepsilon(\alpha,k)\right)$. 
Shatashvili-Vafa's $G_{2}$ superconformal algebra is a quotient of this algebra for $(c,\varepsilon)=(21/2,0)$ modulo an ideal generated in conformal weight $\frac{7}{2}$ (cf.\! Remark \ref{remarkabouttheideal}), in particular, the explicit commutation relations obtained in \cite{Shatashvili-Vafa95} are an artifact of the free field realization the authors used \cite{OFarrill97}. Solving $(c,\varepsilon)=(21/2,0)$ in terms of $\alpha$ and $k$ there are three solutions: $\{\alpha=1,k=-2/3\}$, $\{\alpha=-2,k=-2/3\}$ and $\{\alpha=-1/2,k=1/3\}$. Precisely for this values of $\alpha$ the superalgebra $D(2,1;\alpha)$ is nothing but the superalgebra $osp(4|2)$, then the Shatashvili-Vafa $G_{2}$ superconformal algebra is a quotient of the quantum Hamiltonian reduction of $osp(4|2)$.

\begin{remark}
The existence of the ideal (\ref{ideal}) can be guessed from the fact that the affine vertex algebra $V_{k}(osp(4|2))$ at level $k\in\{-\frac{2}{3},\frac{1}{3}\}$ is not simple, i.e., contains a non-trivial ideal \cite{GorelikKac07}.
\end{remark}

\sloppy

Listed below are the explicit expressions of all the generators of the Shatashvili-Vafa $G_{2}$ superconformal algebra in the case $\{\alpha=1,k=-2/3\}$. Note that we are using the change of basis (\ref{changeofbasis}).

\fussy

\begin{eqnarray}
G&=&\sqrt{\tfrac{3}{2}}f_{1}+\sqrt{\tfrac{3}{2}}f_{2}+\sqrt{\tfrac{3}{2}} f_{3}+\sqrt{\tfrac{3}{2}}:\Phi^{1}h_{1}:+\sqrt{\tfrac{3}{2}}:\Phi^{2}h_{2}:\nonumber\\
&&+\sqrt{\tfrac{3}{2}}:\Phi^{3}h_{3}:-\sqrt{\tfrac{2}{3}} \partial\Phi^{1}-\sqrt{\tfrac{2}{3}} \partial\Phi^{2}-\sqrt{\tfrac{2}{3}} \partial\Phi^{3},\nonumber
\end{eqnarray}

\begin{eqnarray}
L&=&\tfrac{3}{2}f_{12}+\tfrac{3}{2}f_{13}+\tfrac{3}{2}f_{23}-\tfrac{3}{4}:h_{1}h_{1}:-\tfrac{3}{4} :h_{1}h_{2}:-\tfrac{3}{4}:h_{1}h_{3}:-\tfrac{3}{16}:h_{2}h_{2}:\nonumber\\
&&+\tfrac{3}{8}:h_{2}h_{3}-\tfrac{3}{16} :h_{3}h_{3}:-\tfrac{3}{2} :\Phi^{1}f_{2}:-\tfrac{3}{2}:\Phi^{1}f_{3}:+\tfrac{1}{2}:\Phi^{1}\partial\Phi^{2}:\nonumber\\
&&+\tfrac{1}{2} :\Phi^{1}\partial\Phi^{3}:-\tfrac{3}{2}:\Phi^{2}f_{1}:+3:\Phi^{2}f_{3}:-:\Phi^{2}\partial\Phi^{3}:-\tfrac{3}{2} :\Phi^{3}f_{1}:\nonumber\\
&&+3:\Phi^{3}f_{2}:-\tfrac{1}{2}:\partial\Phi^{1}\Phi^{2}:-\tfrac{1}{2} :\partial\Phi^{1}\Phi^{3}:+:\partial\Phi^{2}\Phi^{3}:+\partial h_{1}\nonumber\\
&&+\tfrac{1}{4} \partial h_{2}+\tfrac{1}{4} \partial h_{3},\nonumber
\end{eqnarray}

\begin{eqnarray}
\Phi&=&3 f_{2}-3 f_{3}-6 :\Phi^{1}\Phi^{2}\Phi^{3}:+3 :\Phi^{2}h_{2}:-3:\Phi^{3}h_{3}:-2 \partial\Phi^{2}+2 \partial\Phi^{3},\nonumber
\end{eqnarray}

\begin{eqnarray}
K&=&3 \sqrt{\tfrac{3}{2}} f_{12}-3 \sqrt{\tfrac{3}{2}} f_{13}-\tfrac{3}{2} \sqrt{\tfrac{3}{2}} :h_{1}h_{2}:+\tfrac{3}{2} \sqrt{\tfrac{3}{2}}:h_{1}h_{3}:-\tfrac{3}{4} \sqrt{\tfrac{3}{2}} :h_{2}h_{2}:\nonumber\\
&&+\tfrac{3}{4} \sqrt{\tfrac{3}{2}} :h_{3}h_{3}:-3 \sqrt{\tfrac{3}{2}} :\Phi^{1}f_{2}:+3 \sqrt{\tfrac{3}{2}}:\Phi^{1}f_{3}:+3 \sqrt{\tfrac{3}{2}}:\Phi^{1}\Phi^{2}h_{1}:\nonumber\\
&&-\tfrac{3}{2} \sqrt{\tfrac{3}{2}}:\Phi^{1}\Phi^{2}h_{2}:+\tfrac{3}{2} \sqrt{\tfrac{3}{2}} :\Phi^{1}\Phi^{2}h_{3}:-3 \sqrt{\tfrac{3}{2}}:\Phi^{1}\Phi^{3}h_{1}:\nonumber\\
&&-\tfrac{3}{2} \sqrt{\tfrac{3}{2}}:\Phi^{1}\Phi^{3}h_{2}:+\tfrac{3}{2} \sqrt{\tfrac{3}{2}}:\Phi^{1}\Phi^{3}h_{3}:-3 \sqrt{\tfrac{3}{2}}:\Phi^{2}f_{1}:+3 \sqrt{6}:\Phi^{2}\Phi^{3}h_{1}:\nonumber\\
&&+3 \sqrt{\tfrac{3}{2}} :\Phi^{2}\Phi^{3}h_{2}:+3 \sqrt{\tfrac{3}{2}}:\Phi^{2}\Phi^{3}h_{3}:-2 \sqrt{6}:\Phi^{2}\partial\Phi^{3}:+3 \sqrt{\tfrac{3}{2}} :\Phi^{3}f_{1}:\nonumber\\
&&-\sqrt{6} :\partial\Phi^{1}\Phi^{2}:+\sqrt{6}:\partial\Phi^{1}\Phi^{3}:-2 \sqrt{6} :\partial\Phi^{2}\Phi^{3}:+\sqrt{\tfrac{3}{2}}\partial h_{2}-\sqrt{\tfrac{3}{2}} \partial h_{3},\nonumber
\end{eqnarray}

\begin{eqnarray}
X&=&-3 f_{23}-\tfrac{3}{8}:h_{2}h_{2}:-\tfrac{3}{4}:h_{2}h_{3}-\tfrac{3}{8} :h_{3}h_{3}:-\tfrac{3}{2} :\Phi^{1}\Phi^{2}h_{2}:\nonumber\\
&&-\tfrac{3}{2}:\Phi^{1}\Phi^{2}h_{3}:-\tfrac{3}{2}:\Phi^{1}\Phi^{3}h_{2}:-\tfrac{3}{2}:\Phi^{1}\Phi^{3}h_{3}:-\tfrac{1}{2}:\Phi^{1}\partial\Phi^{2}:\nonumber\\
&&-\tfrac{1}{2}:\Phi^{1}\partial\Phi^{3}:-6 :\Phi^{2}f_{3}:+3:\Phi^{2}\Phi^{3}h_{2}:-3:\Phi^{2}\Phi^{3}h_{3}:+5 :\Phi^{2}\partial\Phi^{3}:\nonumber\\
&&-6 :\Phi^{3}f_{2}:+\tfrac{5}{2}:\partial\Phi^{1}\Phi^{2}:+\tfrac{5}{2} :\partial\Phi^{1}\Phi^{3}:-5:\partial\Phi^{2}\Phi^{3}:+\tfrac{1}{2}\partial h_{2}+\tfrac{1}{2}\partial h_{3},\nonumber
\end{eqnarray}

\small
\begin{eqnarray}
M&=&-3 \sqrt{\tfrac{3}{2}} f_{123}+3 \sqrt{\tfrac{3}{2}}:h_{1}f_{2}:+3 \sqrt{\tfrac{3}{2}} :h_{1}f_{3}:-\tfrac{3}{2} \sqrt{\tfrac{3}{2}}:h_{2}f_{1}:+3 \sqrt{\tfrac{3}{2}} :h_{2}f_{3}:\nonumber\\
&&-\tfrac{3}{2} \sqrt{\tfrac{3}{2}}:h_{3}f_{1}:+3 \sqrt{\tfrac{3}{2}} :h_{3}f_{2}:+3 \sqrt{6} :\Phi^{1}f_{23}:-\tfrac{3}{2} \sqrt{\tfrac{3}{2}}:\Phi^{1}h_{1}h_{2}:\nonumber\\
&&-\tfrac{3}{2} \sqrt{\tfrac{3}{2}}:\Phi^{1}h_{1}h_{3}:-3 \sqrt{\tfrac{3}{2}}:\Phi^{1}\Phi^{2}f_{1}:+3 \sqrt{\tfrac{3}{2}}:\Phi^{1}\Phi^{2}f_{2}:+9 \sqrt{\tfrac{3}{2}}:\Phi^{1}\Phi^{2}f_{3}:\nonumber\\
&&-\sqrt{6}:\Phi^{1}\Phi^{2}\partial\Phi^{3}:-3 \sqrt{\tfrac{3}{2}}:\Phi^{1}\Phi^{3}f_{1}:+9 \sqrt{\tfrac{3}{2}}:\Phi^{1}\Phi^{3}f_{2}:+3 \sqrt{\tfrac{3}{2}}:\Phi^{1}\Phi^{3}f_{3}:\nonumber\\
&&+\sqrt{6}:\Phi^{1}\partial\Phi^{2}\Phi^{2}:+\sqrt{6} :\Phi^{1}\partial\Phi^{2}\Phi^{3}:+\sqrt{6} :\Phi^{1}\partial\Phi^{3}\Phi^{3}:-\tfrac{1}{2} \sqrt{\tfrac{3}{2}}:\Phi^{1}\partial h_{1}:\nonumber\\
&&+3 \sqrt{6}:\Phi^{2}f_{13}:+3 \sqrt{\tfrac{3}{2}} :\Phi^{2}h_{1}h_{2}:+3 \sqrt{\tfrac{3}{2}}:\Phi^{2}h_{2}h_{3}:+3 \sqrt{6}:\Phi^{2}\Phi^{3}f_{2}:\nonumber\\
&&-3 \sqrt{6}:\Phi^{2}\Phi^{3}f_{3}:-2 \sqrt{6}:\Phi^{2}\partial\Phi^{3}\Phi^{3}:-\tfrac{5}{2} \sqrt{\tfrac{3}{2}}:\Phi^{2}\partial h_{2}:+3 \sqrt{6}:\Phi^{3}f_{12}:\nonumber\\
&&+3 \sqrt{\tfrac{3}{2}}:\Phi^{3}h_{1}h_{3}:+3 \sqrt{\tfrac{3}{2}}:\Phi^{3}h_{2}h_{3}:-\tfrac{5}{2} \sqrt{\tfrac{3}{2}}:\Phi^{3}\partial h_{3}:+\tfrac{5}{2} \sqrt{\tfrac{3}{2}}:\partial\Phi^{1}h_{1}:\nonumber\\
&&+\sqrt{\tfrac{3}{2}} :\partial\Phi^{1}h_{2}:+\sqrt{\tfrac{3}{2}}:\partial\Phi^{1}h_{3}:+\sqrt{6} :\partial\Phi^{1}\Phi^{1}\Phi^{2}:+
\sqrt{6} :\partial\Phi^{1}\Phi^{1}\Phi^{3}:\nonumber\\
&&-\sqrt{6}:\partial\Phi^{2}h_{1}:+\tfrac{1}{2} \sqrt{\tfrac{3}{2}}:\partial\Phi^{2}h_{2}:-\sqrt{6}:\partial\Phi^{2}h_{3}:-2 \sqrt{6} :\partial\Phi^{2}\Phi^{2}\Phi^{3}:\nonumber\\
&&-\sqrt{6}:\partial\Phi^{3}h_{1}:-\sqrt{6} :\partial\Phi^{3}h_{2}:+\tfrac{1}{2} \sqrt{\tfrac{3}{2}}:\partial\Phi^{3}h_{3}:-\tfrac{1}{2} \sqrt{\tfrac{3}{2}} \partial f_{1}\nonumber\\
&&-\tfrac{5}{2} \sqrt{\tfrac{3}{2}} \partial f_{2}-\tfrac{5}{2} \sqrt{\tfrac{3}{2}} \partial f_{3}-\sqrt{\tfrac{2}{3}} \partial^{2}\Phi^{1}+\sqrt{\tfrac{2}{3}}\partial^{2}\Phi^{2}+\sqrt{\tfrac{2}{3}} \partial^{2}\Phi^{3}.\nonumber
\end{eqnarray}
\normalsize

The free field realization of $W_{k}(\mathfrak{g},x,f)$ inside $V_{k}(\mathfrak{h})\otimes F(\mathfrak{g}_{1/2})$ is induced by the cannonical homomorphism $\mathfrak{g}_{\le}\rightarrow\mathfrak{g}_{0}$, then we simply obtain the free field realization by removing the terms that contains a current $v\in\mathfrak{g}_{\le}\backslash \mathfrak{g}_{0}$, i.e., the terms containing $f$'s.

For example in the case $\{\alpha=1,k=-2/3\}$ the generators $G$ and $\Phi$ look as:

\begin{equation}
G=\sqrt{\tfrac{3}{2}}:\Phi^{1}h_{1}:+\sqrt{\tfrac{3}{2}}:\Phi^{2}h_{2}:
+\sqrt{\tfrac{3}{2}}:\Phi^{3}h_{3}:-\sqrt{\tfrac{2}{3}} \partial\Phi^{1}-\sqrt{\tfrac{2}{3}} \partial\Phi^{2}-\sqrt{\tfrac{2}{3}} \partial\Phi^{3},\nonumber
\end{equation}

\begin{equation}\label{expressionforPhi}
\Phi=-6 :\Phi^{1}\Phi^{2}\Phi^{3}:+3 :\Phi^{2}h_{2}:-3:\Phi^{3}h_{3}:-2 \partial\Phi^{2}+2 \partial\Phi^{3}.
\end{equation}

Therefore we have proved:

\sloppy

\begin{theorem}
 \label{thm:main-thm}
Let $V_{-2/3}(\mathfrak{h})$ be the affine vertex algebra of level $-2/3$ associated to $\mathfrak{h}$ with bilinear form $A$, and $F(\mathfrak{g}_{1/2})$ the vertex algebra of neutral free fermions as defined above. The vectors $G$ and $\Phi$ given by the expressions above generate the $SW(\frac{3}{2},\frac{3}{2},2)$ vertex algebra with $c=21/2$ and $\varepsilon=0$ inside $V_{-2/3}(\mathfrak{h})\otimes F(\mathfrak{g}_{1/2})$. This vertex algebra is not simple and dividing by the ideal (\ref{ideal}) we obtain the Shatashvili-Vafa $G_{2}$ superconformal algebra.
\end{theorem}

\fussy

\begin{remark}
Note that $V_{-2/3}(\mathfrak{h})\otimes F(\mathfrak{g}_{1/2})$ is isomorphic (by a linear transformation on the generators) to the vertex algebra of three free Bosons and three free Fermions with inner product minus the inverse of Cartan matrix (\ref{Cartanmatrix}) of $D(2,1;1)\simeq osp(4|2)$.
\end{remark}

\begin{remark}
This free field realization was found by Mallwitz \cite{Mallwitz95} using the most general ansatz on three free superfields of conformal weights $\frac{1}{2}$. By obtaining this realization from the quantum Hamiltonian formalism we can find explicitly the screening operators associated with the reduction as follows.
\end{remark}

First we rescale the currents $h\in V_{k}(\mathfrak{h})$ and consider instead $\bar{h}:=\frac{h}{\sqrt{k}}$, therefore $V_{k}(\mathfrak{h})$ is identified as a vertex algebra with the Heisenberg algebra $V_{1}(\mathfrak{h})$ associated to $\mathfrak{h}$.

Let $V_{Q}$ denote the lattice vertex algebra \cite{Kac96} associated to the root lattice $Q$  (that correspond to the Cartan matrix that we have fixed at the beginning of the section) of $D(2,1;\alpha)$, i.e., we have three odd simple roots $\{\alpha_{1},\alpha_{2},\alpha_{3}\}$. Then for every lattice element $\alpha$ we have a $V_{1}(\mathfrak{h})$-module $M_{\alpha}$ and a vertex operator $\Gamma_{\alpha}$ which is an intertwiner of type $\binom{M_{0}}{M_{0} \; M_{\alpha}}$, hence its zero mode maps $V_{1}(\mathfrak{h})=M_{0}\rightarrow M_{\alpha}$. 

 Let $M_{-\alpha_{i}/\sqrt{k}}$ be the $V_{1}(\mathfrak{h})$-module with highest weight $-\alpha_{i}/\sqrt{k}$ and $\Gamma_{-\alpha_{i}/\sqrt{k}}$ the intertwiner constructed just as in the lattice case, so that

\begin{equation}
 \left[{\bar{h_{i}}}_{\lambda}\Gamma_{-\alpha_{j}/\sqrt{k}}\right]=-\frac{(\alpha_{i},\alpha_{j})}{\sqrt{k}}\Gamma_{-\alpha_{j}/\sqrt{k}},\;\;\;\; \partial\left(\Gamma_{-\alpha_{j}/\sqrt{k}}\right)=-\bar{h}_{j}\Gamma_{-\alpha_{j}/\sqrt{k}}.\nonumber 
\end{equation}

Define the operators
\begin{equation}
Q_{i}=:\Phi_{i}\Gamma_{-\alpha_{i}/\sqrt{k}}:\,\in V_{1}(\mathfrak{h})\otimes F(\mathfrak{g}_{1/2})\rightarrow M_{-\alpha_{i}/\sqrt{k}}\otimes F(\mathfrak{g}_{1/2}), \; i=1,2,3.\nonumber
\end{equation}

A straightforward computation using \cite{Thielemans91} shows that

\begin{equation}
W_{k}(\mathfrak{g},x,f)\simeq\bigcap_{i=1}^{3}Ker \;{Q_{i}}_{(0)}\subset V_{1}(h)\otimes F(\mathfrak{g}_{1/2}),\nonumber
\end{equation}

\noindent
equals the free field realization of $W_{k}(\mathfrak{g},x,f)$ inside $V_{k}(\mathfrak{h})\otimes F(\mathfrak{g}_{1/2})$ that we have produced above.

\begin{remark}
In fact a similar result can be obtained for the quantum Hamiltonian reduction of any simple Lie superalgebra when the nilpotent $f$ is \emph{super-principal}, that is, there exists an odd nilpotent $F \in \mathfrak{g}_{-1/2}$ with $[F,F]=f$ ($f \in \mathfrak{g}_{-1}$ being a principal nilpotent)  and these two vectors together with $x$ form part of a copy of $osp(1|2) \subset \mathfrak{g}$. Not all Lie superalgebras admit a superprincipal embedding, in particular, it is necessary to admit a root system with all odd simple roots. In this case, one takes $F = \sum_i e_{-\alpha_i}$ the sum of all simple root vectors.  The list of simple Lie superalgebras admiting an $osp(1|2)$ superprincipal embedding consists of
\begin{gather*}
sl(n \pm 1|n), \qquad osp(2n\pm 1|2n), \qquad osp(2n|2n), \\
osp(2n +2 | 2n), \qquad D(2,1;\alpha) 
\end{gather*}

In these case we see that $\mathfrak{g}_{1/2}$ is naturally isomorphic to $\Pi \mathfrak{h}^*$ and we can form the Boson-Fermion system and the screening charges as above. The intersection of their kernels coincides with the quantum Hamiltonian reduction for generic levels. 
\end{remark}

\appendix 

\section{$\lambda$-brackets of the  Shatashvili-Vafa $G_{2}$ superconformal algebra}\label{App:AppendixA}

\begin{equation}
[{\Phi}_{\lambda}\Phi]=(-\frac{7}{2}) \lambda^{2}+6X,\;\;\; [{\Phi}_{\lambda}X]=-\frac{15}{2}\Phi\lambda-\frac{5}{2}\partial\Phi,\nonumber
\end{equation}

\begin{equation}
[{X}_{\lambda}X]=\frac{35}{24}\lambda^{3}-10X\lambda-5\partial X,\;\;\; [{G}_{\lambda}\Phi]=K,\nonumber
\end{equation}

\begin{equation}
[{G}_{\lambda}X]=-\frac{1}{2}G\lambda+M,\;\;\;[{G}_{\lambda}K]=3\Phi\lambda+\partial\Phi,\nonumber
\end{equation}

\begin{equation}
[{G}_{\lambda}M]=-\frac{7}{12}\lambda^{3}+\left(L+4X\right)\lambda+\partial X,\;\;\; [{\Phi}_{\lambda}K]=-3G\lambda-3\left(M+\frac{1}{2}\partial G\right),\nonumber
\end{equation}

\begin{equation}
[{\Phi}_{\lambda}M]=\frac{9}{2}K\lambda-\left(3:G\Phi:-\frac{5}{2}\partial K\right),\;\;\;[{X}_{\lambda}K]=-3K\lambda+3\left(:G\Phi:-\partial K\right),\nonumber
\end{equation}

\begin{equation}
[{X}_{\lambda}M]=-\frac{9}{4}G\lambda^{2}-\left(5M+\frac{9}{4}\partial G\right)\lambda+\left(4:GX:-\frac{7}{2}\partial M-\frac{3}{4}\partial^{2}G\right),\nonumber
\end{equation}

\begin{equation}
[{K}_{\lambda}K]=-\frac{21}{6}\lambda^{3}+6\left(X-L\right)\lambda+3\partial\left(X-L\right),\nonumber
\end{equation}

\begin{equation}
[{K}_{\lambda}M]=-\frac{15}{2}\Phi\lambda^{2}-\frac{11}{2}\partial\Phi\lambda+3\left(:GK:+2:L\Phi:\right),\nonumber
\end{equation}

\begin{eqnarray}
[{M}_{\lambda}M]&=&-\frac{35}{24}\lambda^{4}+\frac{1}{2}\left(20X-9L\right)\lambda^{2}+\left(10\partial X-\frac{9}{2}\partial L\right)\lambda+\left(\frac{3}{2}\partial^{2}X\right.\nonumber\\
&&\left.-\frac{3}{2}\partial^{2}L-4:GM:+8:LX:\right),\nonumber
\end{eqnarray}

\begin{equation}
[{L}_{\lambda}X]=-\frac{7}{24}\lambda^{3}+2X\lambda+\partial X,\;\;\; [{L}_{\lambda}M]=-\frac{1}{4}G\lambda^{2}+\frac{5}{2}M\lambda+\partial M. \nonumber
\end{equation}

\section{The $SW(\frac{3}{2},\frac{3}{2},2)$ superconformal algebra}\label{App:AppendixB}

Here we follow the presentation in \cite{Noyvert02}. The $SW(\frac{3}{2},\frac{3}{2},2)$ superconformal algebra has six generators $\{G,L,H,\tilde{M},W,U\}$ where $G$ and $L$ generate the $N=1$ superconformal algebra of central charge $c$, and $(H,\tilde{M})$ and $(W,U)$ are two superconformal multiplets of dimensions $\frac{3}{2}$ and $2$ respectively. 

A superconformal multiplet $\hat{\Phi}=\left(\Phi,\Psi\right)$ of dimension $\Delta$ is a pair of two 
primary fields of conformal weights $\Delta$ and $\Delta+\frac{1}{2}$ respectively, such that the $\lambda$-brackets with the supersymmetry
generator $G$ are as follow:

\begin{equation}
[G_{\lambda}\Phi]=\Psi,\;\;\;\;\; [G_{\lambda}\Psi]=(\partial+2\Delta\lambda)\Phi.\nonumber
\end{equation}

The other $\lambda$-brackets between the generators are as follow:

\begin{equation}\label{bracketHwithH}
[H _{\lambda} H]= \frac{c}{3} \lambda^{2}+ \varepsilon \tilde{M}+ 2L +\frac{4}{3}\mu W,
\end{equation}

\begin{equation}
[H _{\lambda} \tilde{M}]=(3G+3\varepsilon H)\lambda +\frac{-2}{3}\mu U+\partial G +\varepsilon\partial H,\nonumber
\end{equation}

\begin{equation}
[\tilde{M} _\lambda \tilde{M}]=\frac{1}{3}c \lambda^3+ (4\varepsilon \tilde{M}+ 8L+ \frac{4}{3}\mu W)\lambda
 +2\varepsilon \partial \tilde{M} +4\partial L + \frac{2}{3}\mu \partial W,\nonumber
\end{equation}

\begin{equation}
[H _\lambda W]=\mu H \lambda +\frac{\varepsilon}{2}U+\frac{\mu}{3}\partial H,\nonumber
\end{equation}

\begin{equation}
[\tilde{M} _\lambda W]=(\frac{\mu}{3}\tilde{M}+2\varepsilon W)\lambda+ \frac{9\mu}{2c}:GH:+\frac{\mu(-27+2c)}{12c}\partial \tilde{M} +\varepsilon \partial W,\nonumber 
\end{equation}

\begin{equation}
[H _\lambda U]=(\frac{-2}{3}\mu \tilde{M} + 2\varepsilon W)\lambda+ \frac{9\mu}{2c}:GH:-\frac{\mu(27+2c)}{12c}\partial \tilde{M} + \frac{\varepsilon}{2}\partial W,\nonumber
\end{equation} 
\small
\begin{equation}
[\tilde{M} _\lambda U]=\mu H \lambda^2+(\frac{5}{2}\varepsilon U + \frac{2}{3}\mu \partial H)\lambda-\frac{9\mu}{2c}:G\tilde{M}:+\frac{9\mu}{c}:LH: +\varepsilon\partial U+ \frac{\mu(-27+2c)}{12c}\partial^2H,\nonumber
\end{equation}
\normalsize

\begin{equation}
[W _\lambda W]=\frac{c}{12}\lambda^3+ (2L+\frac{\varepsilon}{2}\tilde{M}+\frac{\mu(10c-27)}{6c}W)\lambda+ \partial L+ \frac{\varepsilon}{4}\partial \tilde{M} + \frac{\mu(10c-27)}{12c}\partial W,\nonumber
\end{equation}

\begin{eqnarray}
[W _\lambda U]&=&(-\frac{3}{2}G-\frac{3}{4}\varepsilon H)\lambda^2+(\frac{\mu(-27+10c)}{12c}U-\partial G-\frac{\varepsilon}{2}\partial H )\lambda\nonumber\\
&&-\frac{1}{48c} \left(162\varepsilon :G\tilde{M}:+432 \mu :GW:-324:H\tilde{M}:+648:LG:\right.\nonumber\\
&&\left.+324\varepsilon :LH:-8\mu(27+2c)\partial U+6(-27+2c)\partial^2G\right.\nonumber\\
&&\left.+3(-27+2c)\varepsilon \partial^2 H\right),\nonumber
\end{eqnarray}

\begin{eqnarray}
[U _\lambda U]&=&-\frac{c}{12}\lambda^4-(\frac{5}{4}\varepsilon \tilde{M}+5L+ \frac{\mu(-27+10c)}{6c}W)\lambda^2-\left(\frac{5}{4}\varepsilon \partial \tilde{M}+5\partial L\right.\nonumber\\
&&\left.+ \frac{\mu(-27+10c)}{6c}\partial W\right)\lambda-\frac{1}{16c}\left(-144\mu :GU:-108:G\partial G:\right.\nonumber\\
&&\left.-54\varepsilon :G\partial H: + 108 :H\partial H: -108:\tilde{M}\tilde{M}:+216\varepsilon :L\tilde{M}: \right.\nonumber\\
&&\left.+432 :LL:+288\mu :LW:+54 \varepsilon :\partial GH:-3(9-2c)\varepsilon \partial^2\tilde{M}\right.\nonumber\\
&&\left.+24c\partial^2L-4\mu(27-2c)\partial^2W\right),\nonumber
\end{eqnarray}

where  $c,\varepsilon\in \mathbb{C}$ and $\mu=\sqrt{\frac{9c(4+\varepsilon^{2})}{2(27-2c)}}$.

\begin{remark}\label{remarkabouttheideal}
For $(c,\varepsilon)=(\frac{21}{2},0)$  it was checked in \cite{Noyvert02} that $SW(\frac{3}{2},\frac{3}{2},2)$ coincides with the Shatashvili-Vafa $G_{2}$ algebra at central charge $\frac{21}{2}$ modulo the ideal generated by:

\begin{equation}\label{ideal}
2\sqrt{14}:GW:-3:H\tilde{M}:+2:LG:-2\sqrt{14}\partial U.
\end{equation} 

The existence of this ideal was first observed in \cite{OFarrill97}. The relations between the generators of $SW(\frac{3}{2},\frac{3}{2},2)$ in the case $(\frac{21}{2},0)$ and the generators of the Shatashvili-Vafa $G_{2}$ algebra as presented in the Appendix \ref{App:AppendixA} are given by:
\begin{equation}\label{changeofbasis}
\Phi=iH,\;\; K=i\tilde{M},\;\; X=-(L+\sqrt{14}W)/3,\;\; M=-(\partial G+2\sqrt{14}U)/6.
\end{equation}
\end{remark}

\end{document}